\documentclass[12pt]{article}
\usepackage[utf8]{inputenc}
\usepackage{amssymb,amsmath,amsfonts,eucal,mathrsfs,amsthm} 
\usepackage[colorinlistoftodos]{todonotes}
\usepackage{xcolor}
\usepackage{enumitem}
\usepackage{hyperref}
\newtheorem{theorem}{Theorem}
\newtheorem{proposition}[theorem]{Proposition}
\newtheorem{lemma}[theorem]{Lemma}

\newtheorem{corollary}[theorem]{Corollary}

\theoremstyle{definition}
\newcommand{\R}{\mathbb{R}}
\newcommand{\Z}{\mathbb{Z}}

\newcommand{\Sf}{\mathbb{S}}

\newcommand{\Ric}{\mbox{Ric}}

\newcommand{\End}{\mbox{End}}
\newcommand{\Hom}{\mbox{Hom}}
\newcommand{\trace}{\mbox{tr}}

\def\<{{\langle}}
\def\>{{\rangle}}
\def\B{\mathcal{B}}

\def\n{\nabla}

\def\a{\alpha}

\def\be{\begin{equation} }
\def\ee{\end{equation} }

\begin{document}

\title{Ricci pinched compact hypersurfaces\\ in spheres.}
\author{M. Dajczer, M. I. Jimenez and Th. Vlachos}
\date{}
\maketitle

\begin{abstract}
We investigate the topology of the compact hypersurfaces 
in round spheres focusing on those with the Ricci curvature  
satisfying an appropriate bound determined solely by the 
mean curvature of the submanifold. In this paper, the 
application of the Bochner technique yields more robust 
results compared to those presented in \cite{DV} for 
submanifolds that lay in any codimension.
\end{abstract}

Compact submanifolds in the unit sphere $f\colon M^n\to\Sf^{n+p}$, 
$n\geq 4$, with any codimension $p$ have been investigated in 
\cite{DV} under a pinching condition on the Ricci curvature that
depends solely on the norm of the mean curvature vector field. 
In this paper, the special case of codimension $p=1$ is considered. 
Many examples of hypersurfaces meeting the pinching condition 
have been given in \cite{DV}. In this paper, we are able to obtain 
more robust results than in \cite{DV} because the Bochner 
technique applies to hypersurfaces in opposition to the case 
of higher codimension.\vspace{1ex}

Let $f\colon M^n\to\Sf^{n+1}$, $n\geq 4$, be an isometric 
immersion of an $n$-dimensional compact manifold into the 
unit sphere. Given an integer $k$ that satisfies 
$2\leq k\leq n/2$, we denote
$$
b(n,k,H)=
\frac{n(k-1)}{k}+\frac{n(k-1)H}{2k^2}\big(nH+\sqrt{n^2H^2
+4k(n-k)}\big)
$$
being $H$ the length of the (normalized) mean curvature 
vector field of $f$. Throughout the paper, we assume that 
at any point of $M^n$ the (not normalized) Ricci curvature 
satisfies the pinching condition
\be
\Ric_M\geq b(n,k,H)\tag{$\ast$}.
\ee

We say that $(*)$ is \emph{satisfied with equality} at 
$x\in M^n$ if the inequality at that point is not strict, 
that is, there exists a unit vector $X\in T_xM$ such that 
$\Ric_M(X)=b(n,k,H)$. If it happens otherwise, we say that 
$(*)$ is \emph{strict} at $x\in M^n$.
\vspace{1ex}

Recall that the generalized Clifford torus is the standard 
embedding of 
$\mathbb{T}^n_p(r)=\Sf^p(r)\times\Sf^{n-p}(\sqrt{1-r^2})$,
$2\leq p\leq n-2$, into the unit sphere $\Sf^{n+1}$, where 
$\Sf^p(r)$ denotes the $p$-dimensional sphere of radius $r<1$. 
Computations given in \cite{V} yield that $(*)$ is satisfied 
if $p=k$ and $(k-1)/(n-2)\leq r^2\leq k/n$, and that this 
happens with equality.

\begin{theorem}\label{hyp1} Let $f\colon M^n\to\Sf^{n+1}$, 
$n\geq 4$, be an isometric immersion of a compact manifold. 
Assume that $f$ satisfies the pinching condition $(*)$ for 
some $k\geq 2$ where $k<n/2$ if $n$ is even and $k<(n-1)/2$ 
if $n$ is odd. Then $M^n$ is simply connected, hence
orientable, and one of the following cases occurs:
\vspace{1ex}

\noindent $(i)$ The homology groups satisfy 
$$
H_i(M^n;\Z)=H_{n-i}(M^n;\Z)=0\;\,\text{for all}\;\,1\leq i\leq k
$$  
and $H_{n-k-1}(M^n;\Z)=\Z^{\beta_{k+1}(M)}$, where
$\beta_{k+1}(M)$ denotes the $(k+1)$-th Betti number of $M^n$. 
This is necessarily the case if $(*)$ is strict at some point.
\vspace{1ex}

\noindent $(ii)$ The homology groups satisfy 
$$
H_i(M^n;\Z)=H_{n-i}(M^n;\Z)=0\;\,\text{for all}\;\,1\leq i\leq k-1,
$$  
$H_k(M^n;\Z)\neq 0$ is finite and $H_{n-k}(M^n;\Z)=0$. 
For $k=2$, we also assume that $H\neq 0$ at all points.
Then, at any point 
$$
\lambda(n,k,H)
=\frac{1}{2k}\big(nH+\sqrt{n^2H^2+4k(n-k)}\big)
$$
is a principal curvature whose multiplicity $\ell$ satisfies 
$k\leq\ell\leq n-k-1$.  Moreover, equality holds in $(*)$ on 
the principal distribution $T_\lambda=\ker(A-\lambda I)$ 
where $A$ denotes the shape operator of $f$. 
\vspace{1ex}

\noindent $(iii)$ 
$M^n=\mathbb{T}^n_k(r)$ with $(k-1)/(n-2)\leq r^2\leq k/n$
and $f$ is the standard embedding in $\Sf^{n+1}$.
\end{theorem}

Hypersurfaces obtained in part $(ii)$ admit the parametrization
given by Theorem \ref{dupinpar} below on connected components 
of an open dense subset of $M^n$.
\vspace{1ex}

Let $f\colon M^n\to\Sf^{n+1}$, $n\geq 3$, be a hypersurface 
oriented by the unit normal vector field $\xi$ and let $A$ denote
the associated shape operator. Let $\lambda$ be a principal 
curvature of constant multiplicity $\ell$ with $2\leq\ell<n$ 
and let $T_\lambda=\ker(A-\lambda I)$ be the corresponding 
integrable principal distribution. The associated \emph{focal map} 
$f_\lambda\colon M^n\to\Sf^{n+1}$ is given by
$$
f_\lambda=\cos\sigma f+\sin\sigma\xi\;\;\mbox{where}
\;\;\lambda=\cot\sigma.
$$
Let $L$ be the space of leaves 
$M^n/T_\lambda$. It is a standard fact (cf. Theorem~$3.1$ in\cite {CR}) 
that if $M^n$ is complete, then the focal map factors through an 
immersion $g\colon L\to\Sf^{n+1}$ of the $(n-\ell)$-dimensional 
manifold $L^{n-\ell}$. The submanifold $g$ is called the \emph{focal 
submanifold} associated to $\lambda$.

The cases when $k=n/2$ if $n$ is even and $k=(n-1)/2$ if $n$ is
odd have been considered in \cite{DV} for arbitrary codimension. 
In fact, for the first case we have from there the following 
result reiterated here for the sake of completeness.

\begin{theorem} {\em(\cite{DV})} Let $f\colon M^n\to\Sf^{n+1}$, 
$n\geq 4$, be an isometric immersion of a compact manifold of 
even dimension. Assume that
\be\label{strict1}
\Ric_M\geq(n-2)\big(1+H^2+H\sqrt{1+H^2}\big)
\ee
holds at any point of $M^n$. Then one of the following cases occurs:
\vspace{1ex}

\noindent $(i)$ $M^n$ is homeomorphic to $\Sf^n$ and this 
is necessarily the case if at some point of $M^n$
the inequality \eqref{strict1} is strict.\vspace{1ex}

\noindent $(ii)$ The submanifold is the minimal generalized
Clifford torus $\mathbb{T}^n_{n/2}(1/\sqrt{2})$.
\end{theorem}

For hypersurfaces of odd dimension and $k=(n-1)/2$ the 
following result is quite stronger compared to the one 
following from \cite{DV}.

\begin{theorem}\label{hyp3} Let $f\colon M^n\to\Sf^{n+1}$, 
$n\geq 5$, be an isometric immersion of a compact manifold 
of odd dimension. Assume that it holds that
\be\label{ineq2}
\Ric_M\geq\frac{n(n-3)}{n-1}\Big(1+\frac{H}{(n-1)}
\big(nH+\sqrt{n^2H^2+n^2-1}\big)\Big)
\ee
at any point of $M^n$. 
Then one of the following cases occurs:\vspace{1ex}

\noindent $(i)$ $M^n$ is homeomorphic to $\Sf^n$ and 
this is necessarily the case if at some point of $M^n$ 
the inequality \eqref{ineq2} is strict. 
\vspace{1ex}

\noindent $(ii)$ The homology groups $H_i(M^n;\Z)$, 
$1\leq i\leq n-1$, vanish with the exception of 
$H_{(n-1)/2}(M^n;\Z)=\Z_q$ for some $q>1$. For $n=5$ let 
also assume that $H\neq 0$ at any point of $M^n$. 
Then $n=4r+3$ and $\lambda=\lambda(n,k,H)$ is a principal 
curvature with multiplicity $k$ at any point of $M^n$ and 
equality holds in \eqref{ineq2} on the principal 
distribution $T_\lambda$. Moreover, $M^n$ is diffeomorphic 
to the unit normal sphere bundle of the corresponding focal 
submanifold $g\colon L^{k+1}\to\Sf^{n+1}$, being $L^{k+1}$ 
homeomorphic to $\Sf^{k+1}$. 
\vspace{1ex}

\noindent $(iii)$ $M^n=\mathbb{T}^n_{(n-1)/2}(r)$ with 
$(n-3)/2(n-2)\leq r^2\leq(n-1)/2n$ and $f$ is the standard 
embedding in $\Sf^{n+1}$. 
\end{theorem}

From the Theorem \ref{dupinpar} given below it follows that the 
hypersurfaces in  part $(ii)$ admit a global parametrization. 
\vspace{1ex}

The following is a direct consequence of the preceding result.

\begin{corollary}\label{cor} Let $f\colon M^n\to\Sf^{n+1}$, 
$n\geq 5$, be an isometric immersion of a compact manifold 
of odd dimension that satisfies \eqref{ineq2}. 
If $H_{(n-1)/2}(M^n,\Z)$ is torsion free or if 
$n=4r+1$ for $r\geq 2$ then one of the 
following cases occurs:\vspace{1ex}

\noindent $(i)$ $M^n$ is homeomorphic to $\Sf^n$ and 
this is necessarily the case if at some point of $M^n$ 
the inequality \eqref{ineq2} is strict. 
\vspace{1ex}

\noindent $(ii)$  $M^n=\mathbb{T}^n_{(n-1)/2}(r)$ with 
$(n-3)/2(n-2)\leq r^2\leq(n-1)/2n$ and $f$ is the standard 
embedding in $\Sf^{n+1}$. 
\end{corollary}

\section{The Bochner operator}

We start with some algebraic preliminaries inspired by Savo 
\cite{Sv}. Let $V$ be a real $n$-dimensional vector space of 
dimension $n\geq 3$, equipped with a positive definite inner 
product $\<\,,\,\>$. 
We denote by $\mathrm{End}(V)$ the set of self-adjoint 
endomorphisms of $V$ and by $\Lambda^p V^*$, $1\leq p\leq n$, 
the $\binom np$-dimensional real vector space defined 
as the $p$-th exterior power of the dual vector space 
$V^*=\Hom(V,\R)$ of $V$. 
Let $A^{[p]}\in\End(\Lambda^p V^*)$ be given by 
$$
A^{[p]}\omega(v_1,\dots,v_p)
=\sum_{i=1}^p\omega(v_1,\dots,Av_i,\dots,v_p),
$$
where $\omega\in\Lambda^p V^*$ and $v_1,\dots,v_p\in V$.
Then associated to $A\in\End(V)$ there is the endomorphism 
$T_A^{[p]}\in\End(\Lambda^p V^*)$ 
defined by 
$$
T_A^{[p]}=(\trace A)A^{[p]}-A^{[p]}\circ A^{[p]}
$$
which is self-adjoint with respect to the natural inner 
product $\<\,,\,\>$ in $\Lambda^p V^*$.

Let $\{e_1,\dots,e_n\}$ be an orthonormal basis of $V$ 
and let $\{\theta_1,\dots,\theta_n\}$ be the dual basis. 
For every integer $p$ let $\mathcal I_p$ 
be the set of $p$-multi-indices 
$$
\mathcal I_p=\left\{I=(i_1,\dots,i_p)\colon 
1\leq i_1<\cdots<i_p\leq n\right\}.
$$
For each $I=(i_1,\dots,i_p)\in \mathcal I_p$ set
$\theta_I=\theta_{i_1}\wedge\dots\wedge\theta_{i_p}\;\; 
{\text {and}}\;\;e_I=(e_{i_1},\dots,e_{i_p})$.
For any $I,J\in \mathcal I_p$ we have
$$
\theta_I(e_J)=\begin{cases}
\,1&\text{if\, } I=J \\[1mm]
\,0&\text{if \,otherwise.}
\end{cases}
$$

Since $\{\theta_I\colon I\in\mathcal I_p\}$ is 
an orthonormal basis of $\Lambda^p V^*$, 
given $\omega\in\Lambda^p V^*$ we have 
$\omega=\sum_{I\in \mathcal I_p}a_I\theta_I$
where $a_I=\omega(e_I)$.

\begin{lemma}\label{lem-al}
If $A\in\mathrm {End}(V)$ we have for any $I\in \mathcal I_p$,
$1\leq p\leq n$, that
$$
T_A^{[p]}\theta_I=\Big(\trace A\sum_{i\in\mathbf I}\<Ae_i,e_i\>
-\big(\sum_{i\in\mathbf I}\<Ae_i,e_i\>\big)^2\Big)\theta_I,
$$
where $\mathbf I=\{i_1,\dots,i_p\}$ and  $\{e_1,\dots,e_n\}$ 
is an orthonormal basis of $V$ that 
diagonalises $A$. 
\end{lemma}

\proof We have that 
\be\label{Tp}
A^{[p]}\theta_I=\sum_{J\in \mathcal I_p} 
\left(A^{[p]}\theta_I\right)(e_J)\theta_J.
\ee
Then we compute
$\left(A^{[p]}\theta_I\right)(e_J)$ for any $I,J\in\mathcal I_p$.
If $I=(i_1,\dots,i_p),1\leq i_1<\cdots<i_p\leq n$ and 
$J=(j_1,\dots,j_p),1\leq j_1<\cdots<j_p\leq n$, then
\begin{align*}
\left(A^{[p]}\theta_I\right)(e_J)
&=\theta_I(Ae_{j_1},\dots,e_{j_p})+\dots
+\theta_I(e_{j_1},\dots,Ae_{j_p})\\ 
&=\sum_{s=1,s\notin \mathbf J\smallsetminus\{j_1\}}^n
\<Ae_{j_1},e_s\>\theta_I(e_s,e_{j_2},\dots,e_{j_p})\\
&+\sum_{s=1,s\notin \mathbf J\smallsetminus\{j_2\}}^n
\<Ae_{j_2},e_s\>\theta_I(e_{j_1},e_s,e_{j_3},\dots,e_{j_p})\\
&+\dots+\sum_{s=1,s\notin \mathbf J\smallsetminus\{j_p\}}^n
\<Ae_{j_p},e_s\>\theta_I(e_{j_1},\dots,e_{j_{p-1}},e_s).
\end{align*}
Hence
$$
\left(A^{[p]}\theta_I\right)(e_J)
=\sum_{j\in\mathbf J}\<Ae_j, e_j\>\theta_I(e_J),
$$
and thus
$$
\left(A^{[p]}\theta_I\right)(e_J)=
\begin{cases}
\sum_{i\in\mathbf I}\<Ae_i, e_i\>&\text{if\,} J=I \\[1mm]
\,0&\text{if \,otherwise}.
\end{cases}
$$
Then \eqref{Tp} yields
$$
A^{[p]}\theta_I=\sum_{i\in\mathbf I}\<Ae_i,e_i\>\theta_I
$$
for any $I\in \mathcal I_p$. The proof now
follows from the definition of $T_A^{[p]}$.\vspace{2ex}\qed

Let $M^n$ be an orientable Riemannian manifold of dimension $n$. 
For each integer $0\leq p\leq n$, the Hodge-Laplace operator 
acting on differential $p$-forms is defined by 
$$
\Delta=d\delta+\delta d:\Omega^p(M^n)\to\Omega^p(M^n),
$$
where $d$ and $\delta$ are the differential and the co-differential 
operators, respectively. For $p=0$ the Hodge-Laplace operator is 
just the Laplace-Beltrami operator acting on $0$-forms, that is, 
scalar functions. 

A key element in our methodology revolves around the Bochner 
technique, rooted in the Bochner-Weitzenböck formula.
It states that the Laplacian of every $p$-form 
$\omega\in\Omega^p(M^n)$ on a manifold $M^n$ is given by
\be\label{boch.form}
\Delta\omega=\n^*\n\omega+\B^{[p]}\omega,
\ee
where $\n^*\n$ is the so called rough Laplacian or 
connection Laplacian and 
$\B^{[p]}\colon \Omega^p(M^n)\to\Omega^p(M^n)$
is a certain symmetric endomorphism of the bundle of 
$p$-forms called the \emph{Bochner operator}. 

\begin{proposition}\label{nonneg}
Let $f\colon M^n\to\Sf^{n+1}$, $n\geq 4$, be an isometric 
immersion of a compact oriented manifold satisfying the 
inequality $(*)$ for an integer $2\leq k\leq n/2$. Then the 
Bochner operator $\B^{[k]}$ is nonnegative. 
\end{proposition}

\proof Let $\{e_1,\dots,e_n\}$ be a local orthonormal frame 
of the tangent bundle that diagonalizes the shape operator $A$ 
and let $\{\theta_1,\dots,\theta_n\}$ be the dual frame. Set 
$\theta_I=\theta_{i_1}\wedge\dots\wedge\theta_{i_n}$ for any 
$I=(i_1,\dots,i_k)\in \mathcal I_k$. Then for any $k$-form
$\omega=\sum_{I\in \mathcal I_k}a_I\theta_I$ it follows
from Theorem $1$ in \cite{Sv} that 
$$
\<\B^{[k]}\omega,\omega\>=k(n-k)\|\omega\|^2
+\<T_A^{[p]}\omega, \omega\>.
$$
Using Lemma \ref{lem-al}, the above becomes
\begin{align}
\<\B^{[k]}\omega,\omega\>&=k(n-k)\|\omega\|^2&\nonumber\\
&+\sum_{I\in \mathcal I_k}a^2_I
\big(\trace A\sum_{i\in\mathbf I}\<Ae_i,e_i\> 
-\big(\sum_{i\in\mathbf I}\<Ae_i,e_i\>\big)^2 \big).\nonumber
\end{align}
Hence, we have from the Cauchy-Schwarz inequality that
$$
\<\B^{[k]}\omega,\omega\>\geq k(n-k)\|\omega\|^2
+\sum_{I\in \mathcal I_k}a^2_I\sum_{i\in\mathbf I}
\big(\trace A\<Ae_i,e_i\>-k\|Ae_i\|^2 \big).
$$
Using that ${\rm Ric}(X)=n-1+\trace A\<AX,X\>-\|AX\|^2$, the 
above is written as
$$
\<\B^{[k]}\omega,\omega\>\!\geq k(n-k)\|\omega\|^2
+\sum_{I\in \mathcal I_k}a^2_I\sum_{i\in \mathbf I}
\big(k{\rm {Ric}}(e_i)-k(n-1)-(k-1)\trace A\<Ae_i,e_i\>\big).
$$
From Lemma $6$ in \cite{DV} we have that 
$\trace A\<Ae_i,e_i\>\leq nH\lambda(n,k,H)$. Then using the 
assumption on the Ricci curvature we obtain
$$
\<\B^{[k]}\omega,\omega\>\geq\|\omega\|^2k\big(n-k
+kb(n,k,H)-k(n-1)-n(k-1)H\lambda(n,k,H)\big),
$$
and since $kb(n,k,H)=n(k-1)(1+ H\lambda(n,k,H))$ then
the right-hand-side vanishes.\qed

\section{A parametrization}

In this section, our goal is to provide a parametrization for 
the hypersurfaces in spheres with a principal curvature of 
constant multiplicity at least two. This result of 
independent interest in submanifold theory  will be 
applied in one of the forthcoming proofs.
\vspace{1ex}

Let $g\colon L^{n-\ell}\to \Sf^{n+1}\subset\R^{n+2}$,
$2\leq\ell\leq n-2$, be an isometric
immersion into the unit sphere with unit normal bundle
$\Lambda=\{(x,w)\in N_gL\colon\|w\|=1\}$. Then the  projection
$\Pi\colon\Lambda\to L^{n-\ell}$ given by $\Pi(x,w)=x$ is a
submersion whose vertical distribution is $\mathcal V=\ker\Pi_*$.
Let the gradient of $\tau\in C^\infty(L)$ with $0<\tau<\pi/2$  
satisfy $\|\n\tau\|<1$ at any point of $L^{n-\ell}$. Finally, 
let $\Psi\colon\Lambda\to\Sf^{n+1}$ be the map given by
\be\label{param}
\Psi(x,w)=\exp_{g(x)}\delta(x,w)
\ee
where $\delta(x,w)
=-\tau(x)(g_*\n\tau(x)+\sqrt{1-\|\n\tau(x)\|^2}\,w)$
and $\exp$ stands for the exponential map of $\Sf^{n+1}$.
For simplicity, we make use of the same notation for the
corresponding map when composing $\Psi$ with the inclusion
of $\Sf^{n+1}$ into $\R^{n+2}$. Hence, we may write 
$$
\Psi(x,w)=\cos\tau(x)g(x)-\sin\tau(x)\big(g_*\n\tau(x)
+\sqrt{1-\|\n \tau(x)\|^2}\,w\big).
$$
Observe that $\Psi(\Lambda)$ for constant $\tau$ is the 
boundary of the geodesic tube of radius $\tau$ given by
$\left\{\exp_{g(x)}(-\theta w)\colon 0
\leq\theta\leq\tau,\;(x,w)\in\Lambda\right\}$.

\begin{proposition}\label{psi}
Let $M^n\subset\Lambda$ be the open subset of points 
where the map $\Psi$ is regular.
Then the following assertions hold:\vspace{1ex}

\noindent $(i)$  We have that $(x,w)\in M^n$ if
and only if the self adjoint endomorphism
of $T_xL$ given by
\begin{align*}
P(x,w)Y
=&\,\cos\tau(x)(Y-\<Y,\n\tau(x)\>\n\tau(x))
-\sin\tau(x){\rm Hess}\,\tau (x)Y\\
&+\sin\tau(x)\sqrt{1-\|\n\tau(x)\|^2}A^g_wY
\end{align*}
is nonsingular, where $A^g_w$ is the shape operator of $g$.
\vspace{1ex}

\noindent $(ii)$ The Gauss map
$\eta\colon M^n\to\Sf^{n+1}\subset\R^{n+2}$ of $\Psi$ is
given by
$$
\eta(x,w)=\sin\tau(x)g(x)+\cos\tau(x)\big(g_*\n\tau(x)
+\sqrt{1-\|\n\tau(x)\|^2}\,w\big).
$$

\noindent $(iii)$ The shape operator $A$ of $\Psi|_M$
has $\cot\tau$ as a principal curvature with corresponding 
eigenspace $\mathcal{V}$.
\end{proposition}

\proof  At $(x,w)\in M^n$ for $V\in T_{(x,w)}M$
let $c\colon(-\varepsilon,\varepsilon)\to M^n$ be a curve  
of the form $c(t)=(\gamma(t),w(t))$ so that $c(0)=(x,w)$
and $V=c'(0)=(Z,w'(0))$.
A straightforward computation gives
\begin{align}\label{dpsi}
&\Psi_*(x,w)V\nonumber
=g_*P(x,w)Z
-\sin\tau(x)\alpha_g(Z,\n\tau(x))\\
&-\cos\tau(x)\<Z,\n\tau(x)\>\sqrt{1-\|\n\tau(x)\|^2}w
+\sin\tau(x)\frac{\<{\rm Hess}(\tau)(x)Z,\n\tau(x)\>}
{\sqrt{1-\|\n\tau(x)\|^2}}w\nonumber\\
&-\sin\tau(x)\sqrt{1-\|\n\tau(x)\|^2}\frac{\nabla^\perp w}{dt}(0),
\end{align}
where $\a_g\colon TL\times TL\to N_gL$ is the second fundamental
form of $g$, and then part $(i)$ follows.

Since $\eta(x,w)$ is a unit vector tangent to $\Sf^{n+1}$
at $\Psi(x,w)$  we obtain from \eqref{dpsi} that
$\<\Psi_*(x,w)V,\eta(x,w)\>=0$, and this proves part $(ii)$.

A straightforward computation yields
\begin{align*}
\eta_*&(x,w)V
=\,\sin\tau(x)g_*(Z-\<Z,\n\tau(x)\>\n\tau(x))\\
&+\cos\tau(x)g_*({\rm Hess}(\tau)Z-\sqrt{1-\|\n\tau(x)\|^2}A_wZ))\\
&+\cos\tau(x)\big(\alpha_g(Z,\n\tau(x))
+\sqrt{1-\|\n\tau(x)\|^2}\frac{\nabla^\perp w}{dt}(0))\\
&-\sin\tau(x)\<Z,\n\tau(x)\>\sqrt{1-\|\n\tau(x)\|^2}w
-\cos\tau(x)\frac{\<{\rm Hess}(\tau)Z,\n\tau(x)\>}
{\sqrt{1-\|\n\tau(x)\|^2}}w.  
\end{align*}
If $V\in T_{(x,w)}M$ is a vertical vector the above gives
$$
\eta_*(x,w)V
=\cos\tau(x)\sqrt{1-\|\n\tau(x)\|^2}\frac{\nabla^\perp w}{dt}(0).
$$
On the other hand, it follows from \eqref{dpsi} that
$$
\Psi_*(x,w)V
=-\sin\tau(x)\sqrt{1-\|\n\tau(x)\|^2}\frac{\nabla^\perp w}{dt}(0).
$$
Hence each vertical vector $V\in T_{(x,w)} M$ is a principal
vector with $\cot\tau(x)$ as the corresponding principal
curvature.

It remains to prove that there are no other eigenvectors associated 
to $\cot\tau$. We have to show that any solution of
$$
\eta_*(x,w)V=-\cot\tau\Psi_*(x,w)V
$$
with $V=(Z,w'(0))$ as above satisfies $Z=0$.
Taking the normal component to $g$, we obtain that 
$$
(\eta_*(x,w)V)_{N_{g(x)}L}=-\cot\tau(\Psi_*(x,w)V)_{N_{g(x)}L}.
$$
Then a straightforward computation gives that
$\<Z,\nabla\tau\>=0$. Now taking the tangent component, we have
$$
(\eta_*(x,w)V)_{g_*T_xL}=-\cot\tau(\Psi_*(x,w)V)_{g_*T_xL}
$$
and conclude that $Z=0$, as we wished.
\vspace{2ex}\qed

In the sequel, let $f\colon M^n\to\Sf^{n+1}$, $n\geq 4$, be an
orientable hypersurface in the unit sphere with Gauss map $\eta$
and associated shape operator $A$. Let $\lambda>0$ be a principal
curvature of constant multiplicity $2\leq\ell\leq n-2$ and
corresponding principal curvature distribution
$T_\lambda=\ker(A-\lambda I)$.
It is a standard fact that $T_\lambda$ is integrable and that
the leaves are umbilical submanifolds in $\Sf^{n+1}$ along which  
$\lambda$ is constant. In addition, it is well-known that the 
leaves are complete if $M^n$ is complete.  

The associated focal map
$h\colon M^n\to\Sf^{n+1}\subset\R^{n+2}$ to $\lambda$ is
defined as
\be\label{h}
h=\exp_{f}(\sigma\eta)=\cos\sigma f+\sin\sigma\eta,
\ee
where $\sigma\in C^\infty(M)$, $0<\sigma<\pi/2$, is given
by $\cot\sigma=\lambda$. Let $V\subset M^n$ be an open
saturated subset (i.e., $V$ is a union of maximal leaves of
$T_\lambda$) and $L$ the quotient space of leaves of $V$.
Then $L$ is Hausdorff, and hence is an $(n-\ell)$-dimensional 
manifold, if either $V$ is 
the saturation of some cross section to the foliation or all
leaves through points of $V$ are complete. In the following result 
we assume that
$L^{n-\ell}$ is a manifold and let $\pi\colon V\to L^{n-\ell}$ be
the projection. Hence the focal map factors through an
immersion $g\colon L^{n-\ell}\to\Sf^{n+1}$, that is, $h=g\circ\pi$,
which is called the \emph{focal submanifold} associated to $\lambda$.
If $M^n$ is compact, it follows  from the Corollary given in page
$18$ of \cite{P} that also $L^{n-\ell}$ is a compact manifold.

\begin{theorem}\label{dupinpar} Let $f\colon M^n\to\Sf^{n+1}$,
$n\geq 4$, be an orientable hypersurface with a principal curvature
$\lambda$ of constant multiplicity $\ell$ with
$2\leq\ell\leq n-2$.  Then $\tau\in C^\infty(L)$,
$0<\tau<\pi/2$, given by $\sigma=\tau\circ\pi$ satisfies 
$\|\n\tau\|<1$ at any point of $L^{n-\ell}$. 
Let $\Lambda$ be the unit normal bundle of the focal submanifold $g$.
Then there is a local diffeomorphism $j\colon V\to\Lambda$ such 
that $f|_V=\Psi\circ j$ where  $\Psi\colon\Lambda\to\Sf^{n+1}$
is the map given by \eqref{param}. If $M^n$ is compact then 
$L^{n-\ell}$ is compact, $j\colon M^n\to\Lambda$ 
is a covering map and $\Psi$ is a global parametrization.

Conversely, let $g\colon L^{n-\ell}\to\Sf^{n+1}$,
$2\leq\ell\leq n-2$, be a connected submanifold and let
$M^n\subset\Lambda$ be the open subset of points of the
unit normal bundle of $g$
where the map $\Psi$ given by \eqref{param} is
regular. Let $\tau\in C^\infty(L)$, $0<\tau<\pi/2$, be
such that $\|\n\tau\|<1$ at any point of $L^{n-\ell}$.
Then $\Psi\colon M^n\to\Sf^{n+1}$
has a principal curvature $\lambda=\cot\tau\circ\pi$ of
constant multiplicity $\ell$.
\end{theorem}

\proof Since the converse follows from Proposition \ref{psi}
we only argue for the direct part of the statement. 
Let $V\subset M^n$ be such that the quotient $L^{n-\ell}$ is 
a manifold. Along $g\colon L^{n-\ell}\to\Sf^{n+1}\subset\R^{n+2}$
the vector $\eta(x)$ decomposes as
$$
\eta(x)=\<g(\pi(x)),\eta(x)\>g(\pi(x))+g_*(\pi(x))Z+\delta
$$
for $Z\in T_{\pi(x)}L$ and $\delta\in N_gL(\pi(x))$.
Then that $\<g(x),\eta(x)\>=\sin\tau(\pi(x))$ follows from \eqref{h}.
Also  
$$
\<\eta(x),g_*(\pi(x))\tilde X\>=\<Z,\tilde X\>
$$
for any $X\in T_xV$ and $\tilde X=\pi_*(x)X\in T_{\pi(x)}L$.
Moreover, since
\begin{align*}
\<\eta(x),g_*(\pi(x))\tilde X\>
&=\<\eta(x),h_*(x)X\>=\cos\sigma(x)\<\n\sigma,X\>\\
&=\cos\tau(\pi(x))\<\n\tau(\pi(x)),\tilde X\>
\end{align*}
then $Z=\cos\tau(\pi(x))\n\tau(\pi(x))$. It follows
that
$$
\|\delta\|=\cos\tau(\pi(x))\sqrt{1-\|\n\tau(\pi(x))\|^2}.
$$

We argue that $\|\n\tau\|<1$ at any point of $L^{n-\ell}$
and, in particular, we have that $\|\delta\|>0$. On the
contrary, suppose that at some point $\|\n\tau\|=1$
and let $Y\in TM$ be such that $\pi_*Y=\n\tau$. Then
$\eta=\sin\tau\circ\pi\, g+\cos\tau\circ\pi\,g_*\pi_*Y$.
It follows that $Y(\sigma)=1$ and that
$$
(\sin\tau\circ\pi)(g_*\pi_*Y)=(\cos\tau\circ\pi) g\circ\pi-f.
$$
On the other hand, we have from \eqref{h} that
\begin{align*}
(\sin\tau\circ\pi)(g_*\pi_*Y)&=\sin \sigma h_*Y\\
&=\sin\sigma[f_*(\cos\sigma Y-\sin\sigma AY)
+Y(\sigma)(-\sin\sigma f+\cos\sigma\eta)]
\end{align*}
and thus $AY=\cot\sigma Y$, which is a contradiction.
From \eqref{h} there is a unit vector field
$\delta_1\in \Gamma(N_gL)$ such that
$$
f=(\cos\tau\circ\pi) g\circ\pi-(\sin\tau\circ\pi)
\big((g_*\n\tau)\circ\pi
+\sqrt{1-\|(\n\tau)\circ\pi\|^2}\delta_1\big).
$$
Define $j\colon V\to\Lambda$ by
$j(x)=(\pi(x),\delta_1(\pi(x))$. Then the previous equation 
yields $f(x)=\Psi\circ j (x)$. Hence $j$ is a local 
diffeomorphism and $\Psi$ is regular on $j(V)$.
If $M^n$ is compact we take $V=M^n$, and then $L^{n-k}$ is a 
compact manifold, $j$ is a covering map and 
$\Psi(\Lambda)=f(M)$.\qed

\section{The proofs}

In this section, the proofs of the results stated in the 
introduction are given.
\vspace{2ex}

\noindent \emph{Proof of Theorem \ref{hyp1}:} From 
Theorem $1$ in \cite{DV} we have that $M^n$ is simply connected 
and therefore orientable.  Moreover, either we are in case $(i)$ 
or we have that
$$
H_i(M^n;\Z)=0=H_{n-i}(M^n;\Z)\;\,\text{for all}\;\,1\leq i\leq k-1,
$$    
$H_k(M^n;\Z)\neq 0$, $H_{n-k}(M^n;\Z)=\Z^{\beta_k(M)}$  and 
$\lambda(n,k,H)$ is a principal curvature of $f$ 
with multiplicity at least $k$ at any point of $M^n$.
\vspace{1ex}

We need to distinguish two cases:\vspace{1ex}

\noindent\emph{Case I}. Suppose that $\beta_k(M)=0$.  
Poincar\'e duality gives $\beta_{n-k}(M)=0$, and therefore
$$
H_k(M^n;\Z)={\rm {Tor}}(H_k(M^n;\Z))\;\,\mbox{and}\;\,
H_{n-k}(M^n;\Z)={\rm {Tor}}(H_{n-k}(M^n;\Z)).
$$
Since $H_{k-1}(M^n;\Z)=0$ we obtain from the universal coefficient 
theorem for cohomology (cf.\,\cite[pg.\ 244 Corollary 4]{Sp}) that
$$
{\rm {Tor}}(H^k(M^n;\Z))={\rm {Tor}}(H_{k-1}(M^n;\Z))=0.
$$
Poincar\'e duality yields 
$$
{\rm {Tor}}(H_{n-k}(M^n;\Z))={\rm {Tor}}(H^k(M^n;\Z))=0
$$
and hence $H_{n-k}(M^n;\Z)=0$. Thus we are in part $(ii)$
since the remaining of the statement follows from the 
aforementioned Theorem $1$ in \cite{DV}.
\vspace{1ex}

\noindent\emph{Case II}. We argue as in the proof of Case 
$I\!I$ of Theorem \ref{hyp1}. Suppose that $\beta_k(M)>0$. 
Then $M^n$ carries a nontrivial harmonic $k$-form $\omega$. 
Taking the scalar product with $\omega$ on both sides of 
\eqref{boch.form} gives
$$
\|\n\omega\|^2+\<\B^{[p]}\omega,\omega \>
+\frac{1}{2}\,\Delta \|\omega\|^2
=\<\Delta\omega,\omega\>
$$
which is given in \cite{Sv} and in \cite{PL} as Lemma $3.4$.
Proposition \ref{nonneg} yields that the Bochner operator 
$\B^{[k]}$ is nonnegative, and since $\Delta\omega=0$ 
then $\Delta\|\omega\|^2\leq 0$. From the maximum principle 
it follows that $\|\omega\|^2$ is a positive constant. 
Hence $\omega$ is parallel. Thus $M^n$ supports a nontrivial 
parallel $k$-form. Then Theorem $4$ in \cite{Sv} gives that 
$f(M)$ is the torus $\mathbb{T}^n_k(r)$, and from \cite{V} we 
have that the mean curvature is $H=(k-nr^2)/nr\sqrt{1-r^2}$
with $r^2\leq k/n$. Since the Ricci curvature of 
$\mathbb{T}^n_p(r)$ in the principal directions attains the 
values $(k-1)/r^2$ and $(n-k-1)/(1-r^2)$, then the condition 
$(*)$ yields that $r^2\geq(k-1)/(n-2)$.
\vspace{2ex}\qed

\noindent \emph{Proof of Theorem \ref{hyp3}:}
By Theorem \ref{hyp1} in \cite{DV} we have that either
$$
H_i(M^n;\Z)=H_{n-i}(M^n;\Z)=0\;\,\text{for all}\;\,1\leq i\leq k,
$$
or
$$
H_i(M^n;\Z)=H_{n-i}(M^n;\Z)=0\;\,\text{for all}\;\,1\leq i\leq k-1,
$$    
$H_k(M^n;\Z)\neq0$ and $H_{k+1}(M^n;\Z)=\Z^{\beta_k(M)}$. 

In the first case, the manifold $M^n$ is a homology sphere. 
Theorem \ref{hyp1} in \cite{DV} also yields that 
$M^n$ is simply connected. Thus in this case the Hurewicz 
homomorphisms between the homotopy and homology groups are 
isomorphisms, and hence $M^n$ is a homotopy sphere. Since the 
generalized Poincar\'e conjecture holds due to the work of 
Smale and Freedman then $M^n$ is homeomorphic to $\Sf^n$. 

Hereafter we deal with the second case. Part $(ii)$ of Theorem 
\ref{hyp1} in \cite{DV} yields that $\lambda=\lambda(n,k,H)$ is 
a principal curvature with multiplicity $k$ at any point of $M^n$. 
By Theorem $3$ in \cite{DV} we have that either 
the homology of $M^n$ is isomorphic to the one of  
$\Sf^k\times\Sf^{k+1}$, or  the homology of $M^n$ is 
$H_k(M^n,\Z)=\Z_q$ for some $q>1$, satisfies 
$H_0(M^n,\Z)=H_n(M^n,\Z)=\Z$ and is trivial in all other 
cases. By Proposition $7$ in \cite{DV} we have that 
$M^n$ is diffeomorphic to the total space of a sphere bundle 
$\Sf^k\hookrightarrow\mathsf{E}\xrightarrow{p} L^{k+1}$ 
over the quotient manifold $L^{k+1}=M^n/T_\lambda$. 
Moreover, by Theorem $3$ in \cite{DV} we know that $L^{k+1}$ 
is homeomorphic to $\Sf^{k+1}$.\vspace{1ex}

We need to distinguish two cases:
\vspace{1ex}

\noindent\emph{Case I}. Suppose that $H_k(M^n,\Z)=\Z_q$
for some $q>1$. It follows from part $(ii)$ of Theorem 3 in
\cite{DV} that $n=4r+3$. Since  $\lambda=\lambda(n,k,H)$ is
a principal curvature with multiplicity $k$ at any point of 
$M^n$, then Theorem~\ref{dupinpar} gives that $f$ is globally 
a composition $f=\Psi\circ j$, where $j\colon M^n\to\Lambda$ 
is a covering map and $\Lambda$ is the unit normal bundle 
of the compact associated focal submanifold 
$g\colon L^{k+1}\to\Sf^{n+1}$. It is clear that the compact 
hypersurface $\Psi\colon\Lambda\to\Sf^{n+1}$ satisfies $(*)$. 
Theorem \ref{hyp1} in \cite{DV} yields that the manifolds 
$M^n$ and $\Lambda$ are simply connected and therefore 
$j\colon M^n\to\Lambda$ is a diffeomorphism.
\vspace{1ex}

\noindent\emph{Case II}. Suppose that the homology of 
$M^n$ is isomorphic to the  homology of $\Sf^k\times\Sf^{k+1}$. 
Hence $\beta_k(M)=1$. Then $M^n$ carries a nontrivial 
harmonic $k$-form. Proposition~\ref{nonneg} implies that 
the Bochner operator $\B^{[k]}$ is nonnegative. Since 
$H^k(M^n;\R)\neq 0$ then by Proposition \ref{nonneg} every 
harmonic $k$-form is parallel. Thus the manifold $M^n$ supports 
a nontrivial parallel $k$-form. We have from Theorem $4$ in 
\cite{Sv} that $f(M)$ is the torus $\mathbb{T}^n_k(r)$ whose 
mean curvature is  $H=(k-nr^2)/nr\sqrt{1-r^2}$
with $r^2\leq k/n$. Since the Ricci curvature of 
$\mathbb{T}^n_k(r)$ in the principal directions attains the 
values $(k-1)/r^2$ and $k/(1-r^2)$, then condition 
$(*)$ yields $r^2\geq (k-1)/(n-2)$.\vspace{2ex}\qed

Marcos Dajczer is partially supported by the grant 
PID2021-124157NB-I00 funded by 
MCIN/AEI/10.13039/501100011033/ `ERDF A way of making Europe',
Spain, and are also supported by Comunidad Aut\'{o}noma de la Regi\'{o}n
de Murcia, Spain, within the framework of the Regional Programme
in Promotion of the Scientific and Technical Research (Action Plan 2022),
by Fundaci\'{o}n S\'{e}neca, Regional Agency of Science and Technology,
REF, 21899/PI/22. 

Miguel I. Jimenez is supported by FAPESP with 
the grant 2022/05321-9.

\noindent Marcos Dajczer\\
Departamento de Matemáticas\\ 
Universidad de Murcia, Campus de Espinardo\\ 
E-30100 Espinardo, Murcia, Spain\\
e-mail: marcos@impa.br
\bigskip

\noindent Miguel Ibieta Jimenez\\
Instituto de Ciências Matemáticas e de Computação\\
Universidade de São Paulo\\
São Carlos\\
SP 13566-590 -- Brazil\\
e-mail: mibieta@icmc.usp.br
\bigskip

\noindent Theodoros Vlachos\\
University of Ioannina \\
Department of Mathematics\\
Ioannina -- Greece\\
e-mail: tvlachos@uoi.gr

\end{document}